\documentclass{article}
\usepackage{cite}
\usepackage{graphicx}
\usepackage{dcolumn}

\input{tcilatex}
\begin{document}

\title{On the Majorana solution to the Thomas-Fermi equation}
\date{}
\author{Francisco M. Fern\'{a}ndez\thanks{%
fernande@quimica.unlp.edu.ar} \space and Javier Garcia \\
NIFTA, DQT, Sucursal 4, C. C. 16, \\
1900 La Plata, Argentina}
\maketitle

\begin{abstract}
We analyse the solution to the Thomas-Fermi equation discovered by Majorana.
We show that the series for the slope at origin enables one to obtain
results of accuracy far beyond those provided by available methods. We also
estimate the radius of convergence of this series and conjecture that the
singularity closest to origin is a square-root branch point.
\end{abstract}

\section{Introduction}

\label{sec:intro}

There has been great interest in the accurate calculation of the solution to
the nonlinear differential equation that comes from the Thomas-Fermi model
for neutral atoms\cite{F08,AB11,F11,B13,ABF14,PYDG16,PGDY16,PD17,
PMYD17,PRD18,PKDM19,ZB19, RNBMP19}. Several approaches have been applied for
this purpose; for example: Pad\'{e} Hankel method\cite{F08,AB11,F11,ABF14},
fractional order of rational Euler functions\cite{PYDG16}, fractional order
of rational Bessel functions collocation method\cite{PGDY16}, fractional
order of rational Jacobi functions\cite{PMYD17}, rational Chebyshev functions%
\cite{B13}, fractional order of rational Chebyshev functions of the second
kind\cite{PRD18}, a hybrid approach based on the collocation and
Newton-Kantorovich methods plus fractional order of rational Legendre
functions\cite{PKDM19}, Newton iteration with spectral algorithms based on
fractional order of rational Gegenbauer functions\cite{RNBMP19} and rational
Chebyshev series accelerated through coordinate transformations\cite{ZB19}.
We have just mentioned the most accurate results. Other authors have already
obtained less accurate ones and more often than not reported many wrong
digits as shown in the tables of some of the papers just mentioned\cite
{PYDG16, PGDY16, PMYD17,PD17,PRD18,PKDM19}.

The purpose of this paper is the discussion of a semi-analytical solution to
the Thomas-Fermi equation discovered by Majorana in 1928 that remained
unpublished for a long time as revealed in an enlightening pedagogical
article by Esposito\cite{E02}. Although the Majorana solution was mentioned
in some of the papers just quoted\cite
{B13,ABF14,PYDG16,PGDY16,PMYD17,PD17,PRD18,PMYD17} none of those authors
used it for testing their calculations. They did not realize that this
approach enables one to obtain the slope at origin with any desired accuracy
and in fact 100 significant digits have been reported in Wikipedia
(https://en.wikipedia.org/wiki/Thomas\%E2\%80\%93Fermi\_equation\#cite%
\_note-11). As far as we know, this is the most accurate value of the slope
at origin available nowadays.

One may reasonably argue that theoretical results of such an accuracy have
no physical meaning. However, the Thomas-Fermi equation is commonly chosen
as a benchmark for testing algorithms for solving nonlinear differential
equations numerically and, for this reason, accurate results may be useful.
The purpose of this paper is the analysis of the remarkable accuracy of the
Majorana approach to the Thomas-Fermi equation. In section~\ref{sec:Majorana}
we summarize the main equations shown by Esposito\cite{E02}, in section~\ref
{sec:singular point}~we analyze the convergence properties of the Majorana
series numerically and in section~\ref{sec:conclusions} we summarize the
main results and draw conclusions.

\section{The Majorana transformation}

\label{sec:Majorana}

The Thomas-Fermi equation for neutral atoms can be easily reduced to the
dimensionless nonlinear second-order differential equation
\begin{equation}
\varphi ^{\prime \prime }(x)=x^{-1/2}\varphi (x)^{3/2},\;\varphi
(0)=1,\;\varphi (\infty )=0.  \label{eq:T-F}
\end{equation}
In order to solve this equation one has to determine the unknown slope at
origin $\varphi _{0}^{\prime }=\varphi ^{\prime }(0)$ that is consistent
with the boundary conditions. By means of the change of independent and
dependent variables
\begin{equation}
t=144^{-1/6}x^{1/2}\varphi ^{1/6},\;u=-\left( \frac{16}{3}\right)
^{1/3}\varphi ^{-4/3}\varphi ^{\prime },  \label{eq:change_var}
\end{equation}
Majorana derived the nonlinear first-order differential equation\cite{E02}
\begin{equation}
\frac{du(t)}{dt}=8\frac{tu(t)^{2}-1}{1-t^{2}u(t)},\;u(1)=1.  \label{eq:du/dt}
\end{equation}
It can be proved that the slope at origin is given by
\begin{equation}
\varphi _{0}^{\prime }=-\left( \frac{3}{16}\right) ^{1/3}u(0).
\label{eq:slope}
\end{equation}

The solution $u(t)$ can be expanded as
\begin{equation}
u(t)=\sum_{j=0}^{\infty }a_{j}\tau ^{j},\;a_{0}=1,\;\tau =1-t,
\label{eq:u_series}
\end{equation}
where the coefficients $a_{n}$ satisfy the recurrence relation\cite{E02}.
\begin{eqnarray}
a_{m} &=&\frac{1}{2(m+8)-(m+1)a_{1}}\left\{ \sum_{n=1}^{m-2}a_{m-n}\left[
(n+1)a_{n+1}-2(n+4)a_{n}\right. \right.  \nonumber \\
&&\left. +(n+7)a_{n-1}\right] +\left. \left[ m+7-2(m+3)a_{1}\right]
a_{m-1}+(m+6)a_{1}a_{m-2}\right\}  \nonumber \\
&&  \label{eq:rec_rel}
\end{eqnarray}
and $a_{1}$ is a root of $a_{1}^{2}-18a_{1}+8=0$. If $a_{1}=9-\sqrt{73}$
then the coefficients $a_{m}$ are all positive and Esposito\cite{E02}
estimated $a_{n}/a_{n-1}\sim 4/5$ for $n\rightarrow \infty $. On the other
hand, when $a_{1}=9+\sqrt{73}$ the magnitude of the coefficients $\left|
a_{n}\right| $ appears to increase unboundedly and the ratio $\left|
a_{n}/a_{n-1}\right| $ oscillates. Therefore, one expects to obtain
reasonable results only for the first choice that we consider from now on.

It follows from equations (\ref{eq:slope}) and (\ref{eq:u_series}) that one
can obtain approximate solutions to the slope at origin from the sequence of
partial sums
\begin{equation}
\varphi _{0,N}^{\prime }=-\left( \frac{3}{16}\right)
^{1/3}S_{N},\;S_{N}=\sum_{j=0}^{N}a_{j},\;N=1,2,\ldots .
\label{eq:partial_sums}
\end{equation}
Esposito\cite{E02} estimated $\varphi _{0}^{\prime }=-1.588$ (notice that
there is a misprint in his paper) and this result has been cited by several
authors\cite{PYDG16,PGDY16,PMYD17,PD17,PRD18}. However, none of them tried
to obtain accurate results from equations (\ref{eq:rec_rel}) and (\ref
{eq:partial_sums}), except for the accurate value of the slope at origin in
Wikipedia mentioned above.

\begin{figure}[tbp]
\begin{center}
\includegraphics[width=9cm]{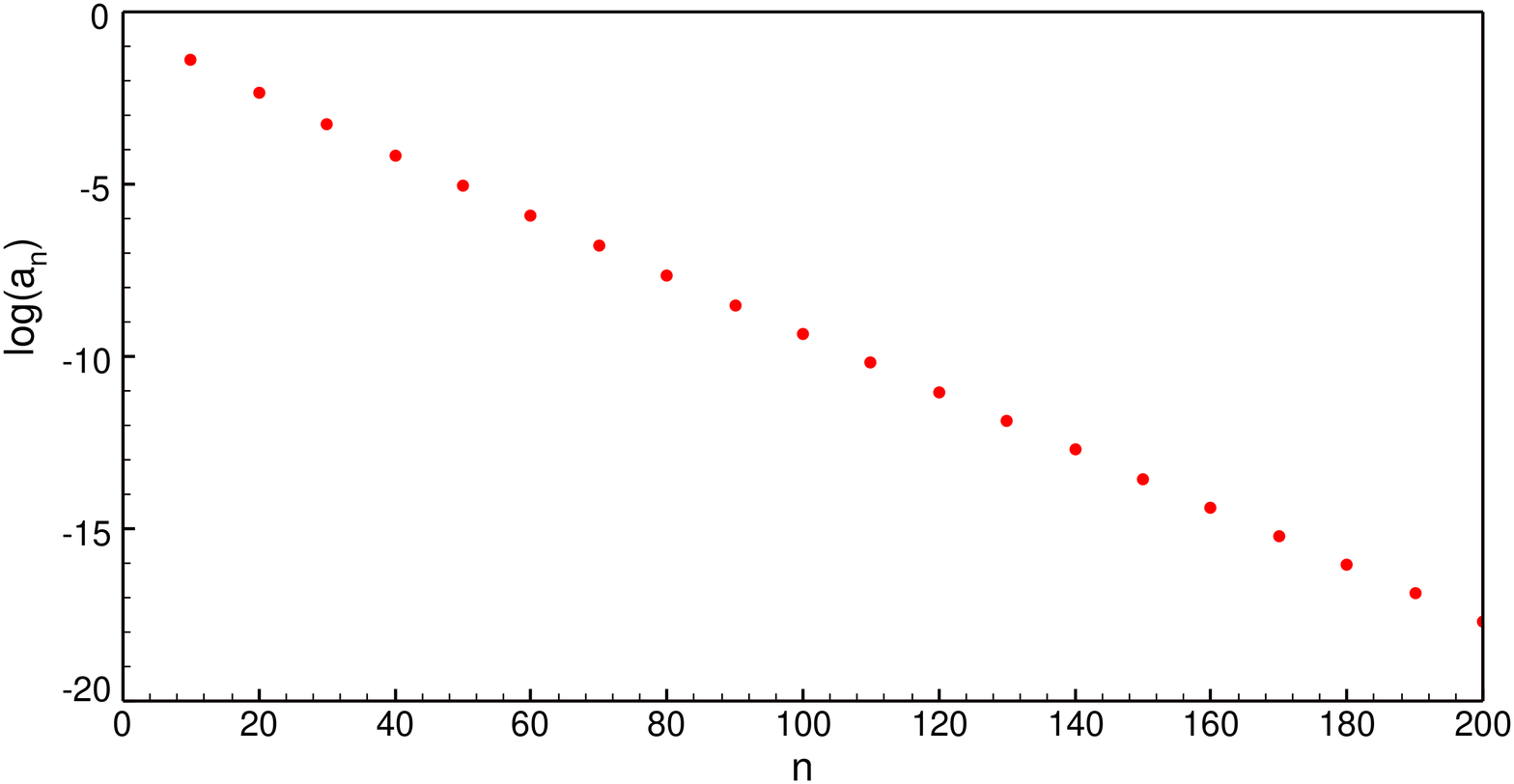}
\end{center}
\caption{Logarithm of the expansion coefficients}
\label{Fig:log(a_n)}
\end{figure}

This approach exhibits two great advantages: first, the expansion
coefficients $a_{n}$ are positive so that $S_{N}<S_{N+1}$ and the sequence
converges from below. Second, the expansion coefficients decrease
exponentially as shown in Figure~\ref{Fig:log(a_n)} and, consequently, the
rate of convergence is remarkable. From the first $N\leq 5000$ partial sums
we estimate
\begin{eqnarray*}
\varphi _{0,5000}^{\prime }
&=&-1.58807102261137531271868450942395010945274662167482 \\
&&561676567741816655196115430926233203397013842866526981 \\
&&402442246654652538163078118812938288132770903659214744 \\
&&086806778901475626169733142185514721772803432126893905 \\
&&784750244483250437472432504160153205506367895247811586 \\
&&772904901159198239550290936872919960125747115454641997 \\
&&589356517736943217510752047768900396683865938312577933 \\
&&839690794052221882937586
\end{eqnarray*}
that is supposed to be accurate to the last digit as we carried out the
calculation numerically with sufficient accuracy (a simple Python program
for this purpose is available at
https://zenodo.org/record/4681779\#.YH2Zorh1Yqg). It is unlikely that any of
the approaches applied to this problem\cite
{F08,AB11,F11,B13,ABF14,PYDG16,PGDY16,PD17,PRD18,PKDM19,ZB19} (and
references therein) can provide a result of such an accuracy. Notice that
the error of a calculation based on $S_{N}$ is roughly of the order of $%
a_{N+1}$.

\section{Analysis of the series}

\label{sec:singular point}

In this section we will try to determine some features of the Majorana
series (\ref{eq:u_series}) numerically. To this end we resort to a method
discussed by Hunter and Guerrieri\cite{HG80} some time ago that we develop
in what follows in a less general and simpler way, more suitable for present
needs.

The function
\begin{equation}
f(x)=A\left( 1-\frac{x}{x_{0}}\right) ^{\nu },  \label{eq:f(x)}
\end{equation}
exhibits a singular point at $x=x_{0}$ for any real exponent $\nu $, except
when it is a positive integer. It is sufficient for present purposes to
consider $x$ real. This function can be expanded in a Taylor series
\begin{equation}
f(x)=\sum_{j=0}^{\infty }f_{j}x^{j},  \label{eq:f(x)_series}
\end{equation}
that converges for all $|x|<\left| x_{0}\right| $. It follows from the
differential equation
\begin{equation}
\left( 1-\frac{x}{x_{0}}\right) f^{\prime }=-\frac{\nu }{x_{0}}f,
\label{eq:f_diff_eq}
\end{equation}
that the expansion coefficients satisfy the recurrence relation
\begin{equation}
(n+1)f_{n+1}=\frac{n}{x_{0}}f_{n}-\frac{\nu }{x_{0}}f_{n},\;n=0,1,\ldots .
\label{eq:f_n_rec_rel}
\end{equation}
We appreciate that the ratio $f_{n+1}/f_{n}$ is a linear function of $%
1/(n+1) $
\begin{equation}
\frac{f_{n+1}}{f_{n}}=\frac{1}{x_{0}}-\frac{\nu +1}{x_{0}(n+1)}.
\label{eq:straight_line_f}
\end{equation}
From equation (\ref{eq:f_n_rec_rel}) we can derive two linear equations
\begin{eqnarray}
x_{0}(n+1)f_{n+1}+\nu f_{n} &=&nf_{n},  \nonumber \\
x_{0}nf_{n}+\nu f_{n-1} &=&(n-1)f_{n-1},  \label{eq:lineqs}
\end{eqnarray}
from which we obtain
\begin{eqnarray}
x_{0} &=&\frac{f_{n}f_{n-1}}{\left( n+1\right) f_{n+1}f_{n-1}-nf_{n}^{2}},
\nonumber \\
\nu &=&\frac{\left( n^{2}-1\right) f_{n+1}f_{n-1}-n^{2}f_{n}^{2}}{\left(
n+1\right) f_{n+1}f_{n-1}-nf_{n}^{2}}.  \label{eq:x_0,nu}
\end{eqnarray}

As stated in section~\ref{sec:Majorana}, the coefficients $a_{n}$ are all
positive and decrease exponentially which suggests that $u(t)$ may exhibit a
singular point as shown in equation (\ref{eq:f(x)}). We can therefore obtain
both the location $\tau _{0}$ and the exponent $\nu $ of the singular point
closest to the origin from the approximate straight line
\begin{equation}
\frac{a_{n+1}}{a_{n}}\approx \frac{1}{\tau _{0}}-\frac{\nu +1}{\tau _{0}(n+1)%
},\;n\gg 1,  \label{eq:straight_line_a}
\end{equation}
as suggested by equation (\ref{eq:straight_line_f}). Alternatively, we may
resort to equation (\ref{eq:x_0,nu}) and estimate these parameters from
\begin{eqnarray}
\tau _{0,n} &=&\frac{a_{n}a_{n-1}}{\left( n+1\right)
a_{n+1}a_{n-1}-na_{n}^{2}},\;\tau _{0}=\lim\limits_{n\rightarrow \infty
}\tau _{0,n},  \nonumber \\
\nu _{n} &=&\frac{\left( n^{2}-1\right) a_{n+1}a_{n-1}-n^{2}a_{n}^{2}}{%
\left( n+1\right) a_{n+1}a_{n-1}-na_{n}^{2}},\;\nu
=\lim\limits_{n\rightarrow \infty }\nu _{n}.  \label{eq:tau_0,nu}
\end{eqnarray}

The second and third columns of Table~\ref{tab:parameters} show values of $%
\tau _{0,n}$ and $\nu _{n}$, respectively, obtained from equation (\ref
{eq:tau_0,nu}) for sufficiently large values of $n$. The convergence is
rather slow but it seems that $\nu _{n}$ converges towards $\nu =1/2$ from
below when $n\rightarrow \infty $. The convergence of these sequences can be
improved by means of of Aitken extrapolation\cite{DB74} and results obtained
from the last $15$ entries in the third column of table~\ref{tab:parameters}
confirms this conjecture. The fourth column shows results for
\begin{equation}
\tau _{0,n}\left( \nu =1/2\right) =\frac{n-1/2}{n+1}\frac{a_{n}}{a_{n+1}}.
\label{eq:tau_(0,n)_nu=1/2}
\end{equation}
It seems that $\tau _{0,n}\left( \nu _{n}\right) $ and $\tau _{0,n}\left(
\nu =1/2\right) $ appear to be monotonously decreasing and increasing,
respectively. On assuming that this behaviour already applies to all $n$ we
conjecture that $1.20168605688<\tau _{0}<1.2016860577$. Straightforward
application of Aitken extrapolation to the last $15$ entries in both columns
enables us to sharpen those bounds to $1.2016860571<\tau _{0}<1.2016860575$.
The same analysis on the third column yields $\nu =0.4998$.

On the other hand, from a numerical fit of the ratio $a_{n+1}/a_{n}$ to a
straight line we estimate $\tau _{0}=1.2016860577$ and $\nu =0.4997$ in
agreement with the results above.

\section{Conclusions}

\label{sec:conclusions}

In this paper we have shown that the Majorana transformation brought to
light by Esposito\cite{E02} in a pedagogical paper enables one to obtain the
slope at origin of the solution to the Thomas-Fermi equation with an
accuracy that has not been achieved with any of the methods proposed so far.
The slope at origin is necessary for the application of any approach because
it is the relevant unknown in the nonlinear differential equation just
mentioned. In addition to it, we estimated the parameters that determine the
singular point closest to the origin of the function $u(t)$ appearing in the
Majorana transformation of the dimensionless Thomas-Fermi equation. We
conjectured that the singularity is a square-root branch point and estimated
its location with reasonable precision by means of lower and upper bounds.

\begin{table}[tbp]
\caption{Parameters for the singular point}
\label{tab:parameters}{\tiny {\
\begin{tabular}{llll}
\multicolumn{1}{c}{$n$} & \multicolumn{1}{c}{$\tau_{0,n}(\nu_n)$} &
\multicolumn{1}{c}{$\nu_n$} & \multicolumn{1}{c}{$\tau_{0,n}(\nu=1/2)$} \\
400000 & 1.20168605769264029435361803588 & 0.499705403241764669703589077642
& 1.20168605680760714618188774385 \\
401000 & 1.20168605769043354866722488742 & 0.499706138707816160983219690953
& 1.20168605680981145947900052075 \\
402000 & 1.20168605768824926820209670303 & 0.49970686850423206121412536936 &
1.20168605681199934431717238544 \\
403000 & 1.2016860576860796865697237324 & 0.499707595193070264721085059588 &
1.20168605681417096352440342843 \\
404000 & 1.20168605768392281024420573709 & 0.499708319421314682402898280176
& 1.20168605681632647794611242793 \\
405000 & 1.20168605768178079115740801082 & 0.499709040445317571614381306288
& 1.20168605681846604642319422634 \\
406000 & 1.20168605767965783307100573377 & 0.499709756820917562940514504979
& 1.20168605682058982586743506764 \\
407000 & 1.2016860576775544722516936779 & 0.499710468331971669275852542178 &
1.20168605682269797125288861779 \\
408000 & 1.20168605767546042850468283397 & 0.499711178432852468843090560012
& 1.20168605682479063564593371776 \\
409000 & 1.20168605767338062590953642012 & 0.499711885437180948142052107338
& 1.20168605682686797024999802803 \\
410000 & 1.20168605767131927381271870206 & 0.499712587885926043832104906147
& 1.20168605682893012444449010861 \\
411000 & 1.20168605766927692214546468536 & 0.499713285557507085244966336909
& 1.20168605683097724578137490504 \\
412000 & 1.20168605766724329181078603638 & 0.499713981941132448940288267412
& 1.20168605683300948001349255953 \\
413000 & 1.20168605766522800138146309189 & 0.499714673721611178261717764848
& 1.20168605683502697114177426129 \\
414000 & 1.20168605766322121156205411946 & 0.499715364255153566375974144531
& 1.20168605683702986142906660283 \\
415000 & 1.20168605766123753508189521398 & 0.49971604848615901584789918849 &
1.20168605683901829144174101888 \\
416000 & 1.2016860576592620129791566382 & 0.499716731547383522068454193635 &
1.20168605684099240002975438654 \\
417000 & 1.20168605765730427494019569306 & 0.499717410088710305095020620403
& 1.20168605684295232439513948465 \\
418000 & 1.20168605765535448415903736306 & 0.499718087499160887322102663134
& 1.2016860568448982000997565051 \\
419000 & 1.20168605765342574745544322327 & 0.499718759200494147400208669829
& 1.20168605684683016110596559594 \\
420000 & 1.20168605765150765389984597273 & 0.499719428790372125524826408218
& 1.2016860568487483397506528092 \\
421000 & 1.20168605764960529812976599862 & 0.499720094468778248392634189577
& 1.20168605685065286682152455643 \\
422000 & 1.20168605764771011681430624946 & 0.49972075921493538246937682624 &
1.2016860568525438715540418048 \\
423000 & 1.20168605764583596292330946038 & 0.499721418145819406112338961443
& 1.20168605685442148167339946943 \\
424000 & 1.20168605764397169184576048614 & 0.499722075152253704050236700056
& 1.20168605685628582337037319071 \\
425000 & 1.20168605764212265017319322498 & 0.499722728329628029434508260725
& 1.20168605685813702137183593605 \\
426000 & 1.20168605764028028246218469624 & 0.49972338068372267220982708216 &
1.20168605685997519893836744777 \\
427000 & 1.20168605763845844575564932541 & 0.499724027284788849437199852618
& 1.20168605686180047790432363325 \\
428000 & 1.20168605763664591943599140173 & 0.499724672088646643944149066554
& 1.20168605686361297865286240567 \\
429000 & 1.20168605763484830282540948745 & 0.499725313083687963509247546151
& 1.20168605686541282018374959084 \\
430000 & 1.20168605763305679932369313611 & 0.499725953390862034718500932713
& 1.20168605686720012011167472364 \\
431000 & 1.20168605763128535363451068823 & 0.499726588003858223074707568887
& 1.20168605686897499470303323157 \\
432000 & 1.20168605762952275020995754139 & 0.499727220914813576068324143944
& 1.20168605687073755885098422327 \\
433000 & 1.20168605762777459527635201156 & 0.499727850091775822559729268675
& 1.20168605687248792614135009801 \\
434000 & 1.2016860576260320960994033668 & 0.49972847868434786817072797909 &
1.20168605687422620884996548424 \\
435000 & 1.20168605762430916779681415112 & 0.499729101651275620448833047275
& 1.20168605687595251797812152143 \\
436000 & 1.20168605762259470597020488796 & 0.499729722982543095775678947105
& 1.20168605687766696322636321741 \\
437000 & 1.20168605762089421604880153902 & 0.49973034066476112550458400528 &
1.20168605687936965305979320755 \\
438000 & 1.20168605761919894526367540425 & 0.49973095786305836731123073029 &
1.20168605688106069470346942218 \\
439000 & 1.2016860576175228147547327723 & 0.499731569488391115336882033527 &
1.20168605688274019417703911361 \\
440000 & 1.20168605761585472407406745998 & 0.499732179566977791786651829719
& 1.2016860568844082562677809745 \\
441000 & 1.20168605761420016645747662319 & 0.499732786072290704444252059124
& 1.20168605688606498459441012161
\end{tabular}
} }
\end{table}


\begin{thebibliography}{99}
\bibitem{F08}  F. M. Fern\'{a}ndez, Comment on: ``Series solution to the
Thomas-Fermi equation'' [Phys. Lett. A 365 (2007) 111], Phys. Lett. A 372
(2008) 5258-5260.

\bibitem{AB11}  S. Abbasbandy and C. Bervillier, Analytic continuation of
Taylor series and the boundary value problems of some nonlinear ordinary
differential equations, Appl. Math. Comput. 218 (2011) 2178-2199.

\bibitem{F11}  F. M. Fern\'{a}ndez, Rational approximation to the
Thomas-Fermi equations, Appl. Math. Comput. 207 (2011) 6433-6436.

\bibitem{B13}  J. P. Boyd, Rational Chebyshev Series for the Thomas-Fermi
Function: Endpoint Singularities and Spectral Methods, J. Comp. Appl. Math.
244 (2013) 90-101.

\bibitem{ABF14}  P. Amore, J. P. Boyd, and F. M. Fern\'{a}ndez, Accurate
calculation of the solutions to the Thomas-Fermi equations, Appl. Math.
Comput. 232 (2014) 929-943.

\bibitem{PYDG16}  K. Parand, H. Yousefi, M. Delkhosh, and A. Ghaderi, A
novel numerical technique to obtain an accurate solution to the Thomas-Fermi
equation, Eur. Phys. J. Plus 131 (2016) 228.

\bibitem{PGDY16}  K. Parand, A. Ghaderi, M. Delkhosh, and H. Yousefi, A new
approach for solving nonlinear Thomas-Fermi equation based on fractional
order rational Bessel functions, Elect. J. Diff. Eq. 2016 (2016) 1-18.

\bibitem{PMYD17}  K. Parand, P. Mazaheri, and M. Delkhosh, Fractional order
of rational Jacobi functions for solving the non-linear singular
Thomas-Fermi equation, Eur. Phys. J. Plus 132 (2017) 77.

\bibitem{PD17}  K. Parand and M. Delkhosh, Accurate solution of the
Thomas-Fermi equation using the fractional order of rational Chebyshev
functions, J. Comp. Appl. Math. 317 (2017) 624-642.

\bibitem{PRD18}  K. Parand, K. Raibei, and M. Delkhosh, An efficient
numerical method for solving nonlinear Thomas-Fermi equation, Acta Univ.
Sapientiae Mathematica 10 (2018) 134-151.

\bibitem{PKDM19}  F. A. Parand, Z. Kalantari, M. Delkhosh, and F.
Mirahmadian, A computationally hybrid Method for Solving a famous physical
problem on an unbounded domain, Commun. Theor. Phys. 71 (2019) 9-15.

\bibitem{RNBMP19}  A. H. Hadian-Rasanan, N. Mehran, A. Bahramnezhad, M. M.
Moayeri, and K. Parand, A comparison between pre-Newton and post-Newton
approaches for solving a physical singular second-order boundary problem in
the semi-infinite interval, arXiv:1909.04066 [math.NA]

\bibitem{ZB19}  X. Zhang and J. P. Boyd, Revisiting the Thomas-Fermi
equation: accelerating rational Chebyshev series through coordinate
transformations, Appl. Num. Math. 135 (2019) 186-205.

\bibitem{E02}  S. Esposito, Majorana solution of the Thomas-Fermi equation,
Am. J. Phys. 70 (2002) 852-856.

\bibitem{HG80}  C. Hunter and B. Guerrieri, Deducing the properties of
singularities of functions from their Taylor series coefficients, SIAM J.
Appl. Math. 39 (1980) 248-263.

\bibitem{DB74}  G. Dahlquist and A Bj\"{o}rck, Numerical Methods,
Prentice-Hall, Englewood Cliffs, 1974.
\end{thebibliography}
\end{document}